\documentclass[12pt,a4paper]{amsart}
\usepackage{amssymb,amsfonts,amsthm,graphicx,color,amsmath}

\theoremstyle{plain}
\newtheorem{theorem}{Theorem}
\newtheorem{lemma}{Lemma}

\newtheorem{definition}{Definition}
\newtheorem{corollary}{Corollary}

\numberwithin{equation}{section} \numberwithin{theorem}{section}
\numberwithin{lemma}{section} \numberwithin{definition}{section}
\numberwithin{corollary}{section}
\numberwithin{proposition}{section} \textheight =24cm
\begin{document}
\title[$n$-space $q$-series identities]{Some $n$-space $q$-binomial theorem extensions and similar identities}
\author{Geoffrey B Campbell}
\address{Mathematical Sciences Institute,
         The Australian National University,
         Canberra, ACT, 0200, Australia}

\email{Geoffrey.Campbell@anu.edu.au}


\begin{abstract}
We give an $n$-space generalized $q$-binomial theorem, and some new $q$ series identities that resemble the traditional $q$ series partition generating functions. These identities enumerate stepping stone weighted vector partitions.
\end{abstract}

\maketitle

\section{Introduction} \label{S:intro}

The literature on $q$ series goes way back to the nineteenth century, starting with Heine \cite{eH1847} and \cite{eH1878}. This generalized the classical hypergeometric series work introduced by Gauss \cite{cG1813}. A beautiful account of these $q$ series identities and their proof as well as their application to integer partitions is given in the excellent book by Andrews \cite{gA1976}. Such $q$ series identities, including the celebrated Rogers-Ramanujan identities, were shown in the 1980s by Baxter \cite{rB1982} to be applicable in the calculation of free energy in statistical mechanics models, such as in his exact solution of \cite{rB1982}.

In his classical account of the theory of partitions, Andrews \cite[chapter 2]{gA1976} shows that many of the time honoured partition identities first given by Euler, Gauss, Heine and Jacobi derive from the $q$-binomial theorem originally given by Cauchy \cite{aC1893},

\begin{theorem}
\textbf{The $q$-binomial theorem.} If $\left|q\right|<1, \left|t\right|<1$, for all complex $a$,
\begin{equation}\label{1.1}
  1+\sum_{k=1}^\infty \frac{(1-a)(1-aq)(1-aq^2)...(1-aq^{k-1})t^k}{(1-q)(1-q^2)(1-q^3)...(1-q^k)} = \prod_{k=0}^\infty \frac{1-atq^k}{1-tq^k}.
\end{equation}
  \end{theorem}

We extend the approach in Andrews \cite{gA1976} which derives many generating functions of integer partition identities from cases of the $q$-binomial theorem. We do so by moving the logic into $n$-space generalizations of the $q$-binomial theorem, as well as applying the methodology to the newer, and not so well-known, \emph{visible point vector identities} that appeared in the 1990s. We therefore, give some new identities that are similar to traditional $q$ series pertaining to \emph{stepping stone} type vector partitions in Euclidean $n$-space regions for the lattice point summations. Along the way we use the standard notation for $q$-shifted factorials,

  \[
    (a;q)_{k}=
    \begin{cases}
      1,                               &\text{$k=0$;}\\
      (1-a)(1-aq)...(1-aq^{k-1}),      &\text{$k=1,2,...$,}
    \end{cases}
  \]

so then $(a;q)_{\infty}$ is understood as the limit $k \rightarrow \infty$, with occasional abbreviated form $(a)_k=(a;q)_k$.

\section{Setting up a higher dimensional approach to $q$-binomial theorem} \label{S:Setting}

We begin with the
\begin{definition}
  Define the function $F_n (t)$ for all complex numbers $a,t$
with $|a|,|t|<1$, and for $n\geq 1$ with $|x_1|,|x_2|,...,|x_n|<1$ by the sequence of functions $F_n(t)$ with
\begin{equation}\label{2.1}
  F_0(-;a,t)= \frac{1-at}{1-t}
\end{equation}
and for all of $\left|x_1\right|,\left|x_2\right|,\left|x_3\right|,...,\left|x_n\right|<1$, by
\begin{equation}\label{2.2a}
  F_n(x_1,x_2,x_3,...,x_n;a,t)= \prod_{\alpha_1,\alpha_2,\alpha_3,...,\alpha_n \geq0}\frac{1-x_1^{\alpha_1} x_2^{\alpha_2} ... x_n^{\alpha_n}at}{1-x_1^{\alpha_1} x_2^{\alpha_2} ... x_n^{\alpha_n}t} \equiv \sum_{k=1}^\infty {_{n}}{A_{k}} t^k
\end{equation}
where $\alpha_1,\alpha_2,\alpha_3,...,\alpha_n$ may be all zero.
\end{definition}

By the way, note that by this definition the right side of (\ref{1.1}) is $F_1(q;a,t)$, so it is evident we are working on a generalised $q$-binomial expression. Next, based on this definition we can assert the

\begin{theorem} \label{thm 2.1} The $n$-space $q$-binomial theorem.
\begin{equation}\label{2.3a}
  {_{n}}{A_{k}}(x_1,x_2,x_3,...,x_n;a)= \det(a_{ij})/k!
\end{equation}
where the determinant is of order $k$, and

  \[
    a_{ij}=
    \begin{cases}
      \frac{1-a^{i-j+1}}{(1-{x_1}^{i-j+1})(1-{x_2}^{i-j+1})(1-{x_3}^{i-j+1})...(1-{x_n}^{i-j+1})},             &\text{$i \geq j$;}\\
      -i,      &\text{$i=j-1$;}\\
      0,      &\text{$otherwise$.}
    \end{cases}
  \]
\end{theorem}

\section{Proof of the $n$-space $q$-binomial theorem.} \label{S:$n$-space proof}

Before giving the full $n$-space proof, we mention that the cases $n=2$ and $n=3$ of theorem \ref{2.1} are not in the literature. Also, for $n=1$, the $q$-binomial theorem itself, there is no proof using determinants in the literature. However, we first prove the $n$-space identity, before then giving a diverse range of examples.

\begin{proof}
  We see first of all that
  \begin{eqnarray}
  \nonumber 
    F_n(x_1,x_2,x_3,...,x_n;a,t) &=& \exp \left( \sum_{\alpha_1,\alpha_2,\alpha_3,...,\alpha_n\geq0} \log\left( \frac{1-x_1^{\alpha_1} x_2^{\alpha_2} ... x_n^{\alpha_n}at}{1-x_1^{\alpha_1} x_2^{\alpha_2} ... x_n^{\alpha_n}t} \right) \right) \\
  \nonumber
      &=& \exp \left( \sum_{k=1}^{\infty} \frac{1-a^k}{\left(1-{x_1}^k\right) \left(1-{x_2}^k\right) \left(1-{x_3}^k\right) \cdots \left(1-{x_n}^k\right)} \frac{t^k}{k}\right).
  \end{eqnarray}
  Since $F_n(x_1,x_2,x_3,...,x_n;a,t)= F_n(t) = \sum_{k=0}^{\infty} {_{n}}{A_{k}} t^k$, differentiating both sides of the equation
  \begin{equation} \nonumber
    \log \left( \sum_{k=0}^{\infty} {_{n}}{A_{k}} t^k \right)
    =  \sum_{k=1}^{\infty} \frac{1-a^k}{\left(1-{x_1}^k\right) \left(1-{x_1}^k\right) \left(1-{x_2}^k\right) \cdots \left(1-{x_n}^k\right)} \frac{t^k}{k}
  \end{equation}
  with respect to $t$, then equating like powers of $t$, leads to the set of simultaneous equations for any positive integer $k$
\begin{equation}   \nonumber 
  1\times{_{n}}{A_{1}} = \frac{1-a}{\left(1-{x_1}\right) \left(1-{x_2}\right) \cdots \left(1-{x_n}\right)},
\end{equation}
\begin{equation}   \nonumber
  2\times{_{n}}{A_{2}} = \frac{1-a^2}{\left(1-{x_1}^2\right) \left(1-{x_2}^2\right) \cdots \left(1-{x_n}^2\right)} + {_{n}}{A_{1}} \frac{1-a}{\left(1-{x_1}\right) \left(1-{x_2}\right) \cdots \left(1-{x_n}\right)},
\end{equation}
\begin{equation}  \nonumber
  3\times{_{n}}{A_{3}} = \frac{1-a^3}{\left(1-{x_1}^3\right) \left(1-{x_2}^3\right) \cdots \left(1-{x_n}^3\right)} + {_{n}}{A_{1}} \frac{1-a^2}{\left(1-{x_1}^2\right) \left(1-{x_2}^2\right) \cdots \left(1-{x_n}^2\right)}
\end{equation}
\begin{equation}  \nonumber
   + {_{n}}{A_{2}} \frac{1-a}{\left(1-{x_1}\right) \left(1-{x_2}\right) \cdots \left(1-{x_n}\right)}, \\
\end{equation}
\begin{equation}  \nonumber
  4\times{_{n}}{A_{4}} = \frac{1-a^4}{\left(1-{x_1}^4\right) \left(1-{x_2}^4\right) \cdots \left(1-{x_n}^4\right)} + {_{n}}{A_{1}} \frac{1-a^3}{\left(1-{x_1}^3\right) \left(1-{x_2}^3\right) \cdots \left(1-{x_n}^3\right)} \\
\end{equation}
\begin{equation}  \nonumber
  + {_{n}}{A_{2}} \frac{1-a^2}{\left(1-{x_1}^2\right) \left(1-{x_2}^2\right) \cdots \left(1-{x_n}^2\right)} + {_{n}}{A_{3}} \frac{1-a}{\left(1-{x_1}\right) \left(1-{x_2}\right) \cdots \left(1-{x_n}\right)}, \\
\end{equation}
   etc., down to the $k$th equation
   \begin{equation}  \nonumber
  k\times{_{n}}{A_{k}} = \frac{1-a^k}{\left(1-{x_1}^k\right) \left(1-{x_2}^k\right) \cdots \left(1-{x_n}^k\right)} + {_{n}}{A_{1}} \frac{1-a^{k-1}}{\left(1-{x_1}^{k-1}\right) \left(1-{x_2}^{k-1}\right) \cdots \left(1-{x_n}^{k-1}\right)} + ... \\
\end{equation}
\begin{equation}  \nonumber
  ... + {_{n}}{A_{k-2}} \frac{1-a^2}{\left(1-{x_1}^2\right) \left(1-{x_2}^2\right) \cdots \left(1-{x_n}^2\right)} + {_{n}}{A_{k-1}} \frac{1-a}{\left(1-{x_1}\right) \left(1-{x_2}\right) \cdots \left(1-{x_n}\right)}. \\
\end{equation}
We note that for the above, the first equation has one unknown, the first and second equations together are \emph{two equations in two unknowns}, the first, second and third equations are \emph{three equations in three unknowns}, and so on. The set of $k$ equations here is a solvable set of \emph{k equations in k unknowns}, and is neatly solved by Cramer's Rule involving determinants. Applying this rule in a straightforward manner yields the theorem.
\end{proof}

The above proof is essentially the same as that given in Macdonald \cite{iM1995} when discussing such combinatorial objects as Schur functions and Jack polynomials. Logarithmic derivatives are taken with respect to $t$, then the recurrence for the coefficients is solved by Cramer's rule.

We write theorem \ref{thm 2.1} as a single equation identity as follows.

\begin{equation}  \label{3.1b}
  F_n(x_1,x_2,x_3,...,x_n;a,t)= \prod_{\alpha_1,\alpha_2,\alpha_3,...,\alpha_n \geq0}\frac{1-x_1^{\alpha_1} x_2^{\alpha_2} ... x_n^{\alpha_n}at}{1-x_1^{\alpha_1} x_2^{\alpha_2} ... x_n^{\alpha_n}t} \\
\end{equation}
\begin{equation}  \nonumber
  = 1 + \frac{1-a}{\left(1-{x_1}\right) \left(1-{x_2}\right) \cdots \left(1-{x_n}\right)} \frac{t}{1!}
\end{equation}
\begin{equation}  \nonumber
+ \begin{vmatrix}
    \frac{1-a}{\left(1-{x_1}\right) \left(1-{x_2}\right) \cdots \left(1-{x_n}\right)} & -1 \\
    \frac{1-a^2}{\left(1-{x_1}^2\right) \left(1-{x_2}^2\right) \cdots \left(1-{x_n}^2\right)} & \frac{1-a}{\left(1-{x_1}\right) \left(1-{x_2}\right) \cdots \left(1-{x_n}\right)} \\
  \end{vmatrix} \frac{t^2}{2!}
  \end{equation}
\begin{equation}  \nonumber
+ \begin{vmatrix}
    \frac{1-a}{\left(1-{x_1}\right) \left(1-{x_2}\right) \cdots \left(1-{x_n}\right)} & -1 & 0 \\
    \frac{1-a^2}{\left(1-{x_1}^2\right) \left(1-{x_2}^2\right) \cdots \left(1-{x_n}^2\right)} & \frac{1-a}{\left(1-{x_1}\right) \left(1-{x_2}\right) \cdots \left(1-{x_n}\right)} & -2 \\
    \frac{1-a^3}{\left(1-{x_1}^3\right) \left(1-{x_2}^3\right) \cdots \left(1-{x_n}^3\right)} & \frac{1-a^2}{\left(1-{x_1}^2\right) \left(1-{x_2}^2\right) \cdots \left(1-{x_n}^2\right)} & \frac{1-a}{\left(1-{x_1}\right) \left(1-{x_2}\right) \cdots \left(1-{x_n}\right)} \\
  \end{vmatrix} \frac{t^3}{3!}
  \end{equation}
\begin{equation}  \nonumber
+ \begin{vmatrix}
    \frac{1-a}{\left(1-{x_1}\right) \left(1-{x_2}\right) \cdots \left(1-{x_n}\right)} & -1 & 0 & 0 \\
    \frac{1-a^2}{\left(1-{x_1}^2\right) \left(1-{x_2}^2\right) \cdots \left(1-{x_n}^2\right)} & \frac{1-a}{\left(1-{x_1}\right) \left(1-{x_2}\right) \cdots \left(1-{x_n}\right)} & -2 & 0 \\
    \frac{1-a^3}{\left(1-{x_1}^3\right) \left(1-{x_2}^3\right) \cdots \left(1-{x_n}^3\right)} & \frac{1-a^2}{\left(1-{x_1}^2\right) \left(1-{x_2}^2\right) \cdots \left(1-{x_r}^2\right)} & \frac{1-a}{\left(1-{x_1}\right) \left(1-{x_2}\right) \cdots \left(1-{x_n}\right)} & -3 \\
    \frac{1-a^4}{\left(1-{x_1}^4\right) \left(1-{x_2}^4\right) \cdots \left(1-{x_n}^4\right)} & \frac{1-a^3}{\left(1-{x_1}^3\right) \left(1-{x_2}^3\right) \cdots \left(1-{x_n}^3\right)} & \frac{1-a^2}{\left(1-{x_1}^2\right) \left(1-{x_2}^2\right) \cdots \left(1-{x_n}^2\right)} & \frac{1-a}{\left(1-{x_1}\right) \left(1-{x_2}\right) \cdots \left(1-{x_n}\right)} \\
  \end{vmatrix} \frac{t^4}{4!}
  \end{equation}
\begin{equation}  \nonumber
 + etc.
\end{equation}

The $1$-space example of equation (\ref{3.1b}) is the $q$-binomial theorem itself, and the following form of this is proven by applying the famous Faà de Bruno's formula (see Abramowitz and Stegun  \cite[chapter 24, pp824]{mA1972}). This form of the $q$-binomial theorem seems to not be in the literature.
\begin{corollary}\label{3.2b}
\begin{equation}\label{3.2}
  F_1(q;a,t)= \prod_{k=0}^\infty \frac{1-atq^k}{1-tq^k}
  = 1 + \frac{1-a}{1-q} \frac{t}{1!}
+ \begin{vmatrix}
    \frac{1-a}{1-q} & -1 \\
    \frac{1-a^2}{1-q^2} & \frac{1-a}{1-q} \\
  \end{vmatrix} \frac{t^2}{2!}
    \end{equation}
  \begin{equation}  \nonumber
+ \begin{vmatrix}
    \frac{1-a}{1-q} & -1 & 0 \\
    \frac{1-a^2}{1-q^2} & \frac{1-a}{1-q} & -2 \\
    \frac{1-a^3}{1-q^3} & \frac{1-a^2}{1-q^2} & \frac{1-a}{1-q} \\
  \end{vmatrix} \frac{t^3}{3!}
+ \begin{vmatrix}
    \frac{1-a}{1-q} & -1 & 0 & 0 \\
    \frac{1-a^2}{1-q^2} & \frac{1-a}{1-q} & -2 & 0 \\
    \frac{1-a^3}{1-q^3} & \frac{1-a^2}{1-q^2} & \frac{1-a}{1-q} & -3 \\
    \frac{1-a^4}{1-q^4} & \frac{1-a^3}{1-q^3} & \frac{1-a^2}{1-q^2} & \frac{1-a}{1-q} \\
  \end{vmatrix} \frac{t^4}{4!}
+ etc.
\end{equation}
  \end{corollary}
    An interesting related aside is that in Mathematica or WolframAlpha, the code

\begin{equation}  \nonumber
Det[\{(1-a)/(1-q), -1,0\}, \\
\{(1-a^2)/(1-q^2), (1-a)/(1-q),-2\},  \\
\end{equation}
\begin{equation}  \nonumber
\{(1-a^3)/(1-q^3), (1-a^2)/(1-q^2), (1-a)/(1-q)\}]  \\
\end{equation}
\begin{equation}  \nonumber
= 6 ((1-a)(1-aq)(1-aq^2))/((1-q)(1-q^2)(1-q^3))
\end{equation}

yields the particular verification

\begin{equation}  \nonumber
\begin{vmatrix}
    \frac{1-a}{1-q} & -1 & 0 \\
    \frac{1-a^2}{1-q^2} & \frac{1-a}{1-q} & -2 \\
    \frac{1-a^3}{1-q^3} & \frac{1-a^2}{1-q^2} & \frac{1-a}{1-q} \\
  \end{vmatrix} =  6\times \frac{(1-a)(1-aq)(1-aq^2)}{(1-q)(1-q^2)(1-q^3)}, \\
  \end{equation}

and the direct application of Faà De Bruno's formula to the general determinant coefficient in equation (3.2) will complete this new proof of $q$-binomial theorem.

\begin{corollary}\label{3.3b}
The 2-space example of equation (\ref{3.1b}) is
\begin{equation}\label{3.3}
  F_2(x,y;a,t)= \prod_{j,k\geq0}\frac{1-x^j y^k at}{1-x^j y^k t} \\
\end{equation}
\begin{equation}  \nonumber
  = 1 + \frac{1-a}{\left(1-x\right) \left(1-y\right)} \frac{t}{1!}
+ \begin{vmatrix}
    \frac{1-a}{\left(1-x\right) \left(1-y\right)} & -1 \\
    \frac{1-a^2}{\left(1-x^2\right) \left(1-y^2\right)} & \frac{1-a}{\left(1-x\right) \left(1-y\right)} \\
  \end{vmatrix} \frac{t^2}{2!}
  \end{equation}
\begin{equation}  \nonumber
+ \begin{vmatrix}
    \frac{1-a}{\left(1-x\right) \left(1-y\right)} & -1 & 0 \\
    \frac{1-a^2}{\left(1-x^2\right) \left(1-y^2\right)} & \frac{1-a}{\left(1-x\right) \left(1-y\right)} & -2 \\
    \frac{1-a^3}{\left(1-x^3\right) \left(1-y^3\right)} & \frac{1-a^2}{\left(1-x^2\right) \left(1-y^2\right)} & \frac{1-a}{\left(1-x\right) \left(1-y\right)} \\
  \end{vmatrix} \frac{t^3}{3!}
  \end{equation}
\begin{equation}  \nonumber
+ \begin{vmatrix}
    \frac{1-a}{\left(1-x\right) \left(1-y\right)} & -1 & 0 & 0 \\
    \frac{1-a^2}{\left(1-x^2\right) \left(1-y^2\right)} & \frac{1-a}{\left(1-x\right) \left(1-y\right)} & -2 & 0 \\
    \frac{1-a^3}{\left(1-x^3\right) \left(1-y^3\right)} & \frac{1-a^2}{\left(1-x^2\right) \left(1-y^2\right)} & \frac{1-a}{\left(1-x\right) \left(1-y\right)} & -3 \\
    \frac{1-a^4}{\left(1-x^4\right) \left(1-y^4\right)} & \frac{1-a^3}{\left(1-x^3\right) \left(1-y^3\right)} & \frac{1-a^2}{\left(1-x^2\right) \left(1-y^2\right)} & \frac{1-a}{\left(1-x\right) \left(1-y\right)} \\
  \end{vmatrix} \frac{t^4}{4!}
  \end{equation}
\begin{equation}  \nonumber
+ etc.
\end{equation}
  \end{corollary}

  The function $F_2(x,y;a,t)$ is interesting for a variety of reasons.  Firstly, $F_2(0,q;a,t)$ and $F_2(q,0;a,t)$ are the same as the right side of equation (\ref{1.1}), and therefore corollary 3.1 is an extension of the $q$-binomial formula.  Secondly, particular cases such as for $|x|<1$

  \begin{equation}  \nonumber
F_2(x,x;a,t)=\prod_{k\geq0} \left(\frac{1-x^kat}{1-x^kt}\right)^k
  \end{equation}

  lead easily to functions like

    \begin{equation}  \nonumber
F_2(q,q;0,qa)=\prod_{k\geq0} \left(\frac{1}{1-q^ka}\right)^k
  \end{equation}

  whose coefficient of  $a^mq^n$   is the number of plane partitions of  $n$  for which the sum of the diagonal parts is  $m$.  (see Andrews \cite[pages 189 and  199]{gA1976}.  Further esoteric restrictions on (\ref{3.2}) yield the celebrated MacMahon generating functions (see Andrews  \cite[page 184]{gA1976}),

   \begin{equation}  \nonumber
F_2(q,q;q^k,q)=\prod_{j\geq1} \left(\frac{1}{1-q^j}\right)^{\min(k,j)},
   \end{equation}

   \begin{equation}  \nonumber
F_2(q,q;0,q)=\prod_{j\geq1} \left(\frac{1}{1-q^j}\right)^j,
   \end{equation}

whose coefficients of  $q^j$   enumerate respectively:-

\begin{enumerate}
  \item the number of  \emph{$k$ rowed} partitions of  $j$,
  \item the number of \emph{unlimited row} plane partitions of  $j$.
\end{enumerate}

$F_2(x,y;a,t)$ can also lead to formulae for the functions  (for $|x|, |y|<1$),

    \begin{equation}  \nonumber
\prod_{j,k\geq0;j+k\neq0} \left(\frac{1}{1-x^jy^k}\right),
   \end{equation}

   \begin{equation}  \nonumber
\prod_{j,k\geq0} \left(1+x^jy^k\right),
   \end{equation}

whose coefficients of  $x^my^n$ generate two dimensional first quadrant Euclidean space vector partitions. Citing the theory in Andrews  \cite[chapter 12]{gA1976})
these coefficients enumerate respectively, the number of partitions of:-
\begin{enumerate}
  \item  \emph{unrestricted} two dimensional vectors summing to vector $\langle m,n \rangle$,
  \item  \emph{distinct} two dimensional vectors summing to vector $\langle m,n \rangle$.
\end{enumerate}

In the present paper we envision each integer lattice point in the vector space as if it were a stepping stone, as the theory for these vector partitions, at least in 2D and 3D Euclidean space, seems to lend itself to this analogy.

  \begin{corollary}\label{3.4a}
The 3-space example of equation (\ref{3.1b}) is
\begin{equation}\label{3.4}
  F_3(x,y,z;a,t)= \prod_{j,k,h\geq0}\frac{1-x^j y^k z^h at}{1-x^j y^k z^h t} \\
\end{equation}
\begin{equation}  \nonumber
  = 1 + \frac{1-a}{\left(1-x\right) \left(1-y\right) \left(1-z\right)} \frac{t}{1!}
+ \begin{vmatrix}
    \frac{1-a}{\left(1-x\right) \left(1-y\right) \left(1-z\right)} & -1 \\
    \frac{1-a^2}{\left(1-x^2\right) \left(1-y^2\right) \left(1-z^2\right)} & \frac{1-a}{\left(1-x\right) \left(1-y\right) \left(1-z\right)} \\
  \end{vmatrix} \frac{t^2}{2!}
  \end{equation}
\begin{equation}  \nonumber
+ \begin{vmatrix}
    \frac{1-a}{\left(1-x\right) \left(1-y\right) \left(1-z\right)} & -1 & 0 \\
    \frac{1-a^2}{\left(1-x^2\right) \left(1-y^2\right) \left(1-z^2\right)} & \frac{1-a}{\left(1-x\right) \left(1-y\right) \left(1-z\right)} & -2 \\
    \frac{1-a^3}{\left(1-x^3\right) \left(1-y^3\right) \left(1-z^3\right)} & \frac{1-a^2}{\left(1-x^2\right) \left(1-y^2\right) \left(1-z^2\right)} & \frac{1-a}{\left(1-x\right) \left(1-y\right) \left(1-z\right)} \\
  \end{vmatrix} \frac{t^3}{3!}
  \end{equation}
\begin{equation}  \nonumber
+ \begin{vmatrix}
    \frac{1-a}{\left(1-x\right) \left(1-y\right) \left(1-z\right)} & -1 & 0 & 0 \\
    \frac{1-a^2}{\left(1-x^2\right) \left(1-y^2\right) \left(1-z^2\right)} & \frac{1-a}{\left(1-x\right) \left(1-y\right) \left(1-z\right)} & -2 & 0 \\
    \frac{1-a^3}{\left(1-x^3\right) \left(1-y^3\right) \left(1-z^3\right)} & \frac{1-a^2}{\left(1-x^2\right) \left(1-y^2\right) \left(1-z^2\right)} & \frac{1-a}{\left(1-x\right) \left(1-y\right) \left(1-z\right)} & -3 \\
    \frac{1-a^4}{\left(1-x^4\right) \left(1-y^4\right) \left(1-z^4\right)} & \frac{1-a^3}{\left(1-x^3\right) \left(1-y^3\right) \left(1-z^3\right)} & \frac{1-a^2}{\left(1-x^2\right) \left(1-y^2\right) \left(1-z^2\right)} & \frac{1-a}{\left(1-x\right) \left(1-y\right) \left(1-z\right)} \\
  \end{vmatrix} \frac{t^4}{4!}
  \end{equation}
\begin{equation}  \nonumber
+ etc.
\end{equation}
  \end{corollary}

The identities in the present paper follow the approach taken by Andrews \cite{gA1976}, who develops his theory of integer partitions starting from an historical viewpoint and uses the $q$-binomial theorem to derive the Euler identities for the number of partitions of a positive integer. He then also uses (\ref{1.1}) to yield the Jacobi triple product, the Heine $q$-series summation and transform, and other $q$-series sum-product identities relevant to partitions.

\section{Example variations for identities generating weighted vector partitions} \label{S:eg variations}

The methods of the book by Macdonald \cite{iM1995} would have determined these formulae
earlier but in his examples he considers only single q-variable and bivariate infinite products.
Although Cramer’s rule (see Birkhoff and Maclaine \cite{gB1977}) is used in examples of this book, the
method has not been specifically applied to products such as those in this paper.

In relation to some researches of the past 53 years, the interesting papers on vector
partition generating functions such as Cheema \cite{mC1966}, Cheema and Motzkin \cite{mC1971}, Gordon \cite{bG1963}
and Wright \cite{eW1956} show that our current methods have not been used in many of the old
problems in vector partition theory.

We continue with a vector partition generating function variation on the work so far, by stating
\begin{definition}
  We define the function $G_n (t)$ for all non-negative integers $n$, and for suitable complex numbers $a$ and
$t$, by
\begin{equation}\label{4.1}
  G_0(-;a,t)= \frac{1-at}{1-t}
\end{equation}
and for all of $\left|x_1\right|,\left|x_2\right|,...,\left|x_n\right|<1$, by
\begin{equation}\label{4.2a}
  G_n(x_1,x_2,...,x_n;a,t,s)= \prod_{\alpha_1,\alpha_2,...,\alpha_n \geq0} \left(\frac{1-x_1^{\alpha_1} x_2^{\alpha_2} ... x_n^{\alpha_n}at}{1-x_1^{\alpha_1} x_2^{\alpha_2} ... x_n^{\alpha_n}t}\right)^{\frac{1}{{\alpha_1}! {\alpha_2}! ... {\alpha_n}!}}
  \end{equation}
  \begin{equation}  \nonumber
  \equiv \sum_{k=1}^\infty {_{n}}{B_{k}} t^k
\end{equation}
where $\alpha_1,\alpha_2,...,\alpha_n$ may be all zero which in that particular case is $G_0(-;a,t)$.
\end{definition}

By the way, note that by this definition the right side $G_1(q;a,t)$, is a variant on the standard $q$-binomial expression (1.1). Next, based on this definition we can assert the

\begin{theorem} \label{thm 4.1} An $n$-space variation on the extended $q$-binomial theorem.
\begin{equation}\label{4.3a}
  {_{n}}{B_{k}}(x_1,x_2,...,x_n;a)= \det(a_{ij})/k!
\end{equation}
where the determinant is of order $k$, and

  \[
    a_{ij}=
    \begin{cases}
      \left(1-a^{i-j+1}\right)\exp\left({x_1}^{i-j+1}+{x_2}^{i-j+1}+...+{x_n}^{i-j+1}\right),             &\text{$i \geq j$;}\\
      -i,      &\text{$i=j-1$;}\\
      0,      &\text{$otherwise$.}
    \end{cases}
  \]
\end{theorem}

\section{Proof of the $n$-space variation of the extended $q$-binomial theorem.} \label{S:variation proof}

As with the full $n$-space $q$-binomial theorem proof, we mention that the cases $n=2$ and $n=3$ of theorem \ref{4.1} are new to the literature. Also, for $n=1$, the $q$-binomial theorem variation, there is no proof in the literature. However, we first prove the $n$-space identity, before then giving a diverse range of examples.

\begin{proof}
  We see that, similar to the earlier proof,
  \begin{eqnarray}
  \nonumber 
    G_n(x_1,x_2,...,x_n;a,t) &=& \exp \left( \sum_{\alpha_1,\alpha_2,...,\alpha_n\geq0} \frac{1}{{\alpha_1}! {\alpha_2}! ... {\alpha_n}!} \log\left( \frac{1-x_1^{\alpha_1} x_2^{\alpha_2} ... x_n^{\alpha_n}at}{1-x_1^{\alpha_1} x_2^{\alpha_2} ... x_n^{\alpha_n}t} \right) \right) \\
  \nonumber
      &=& \exp \left( \sum_{k=1}^{\infty} \left(1-a^{k}\right)\exp\left({x_1}^{k}+{x_2}^{k}+...+{x_n}^{k}\right) \frac{t^k}{k}\right).
  \end{eqnarray}
  Since $G_n(x_1,x_2,x_3,...,x_n;a,t)= G_n(t) = \sum_{k=0}^{\infty} {_{n}}{B_{k}} t^k$, differentiating both sides of the equation
  \begin{equation} \nonumber
    \log \left( \sum_{k=0}^{\infty} {_{n}}{B_{k}} t^k \right)
    =  \sum_{k=1}^{\infty} \left(1-a^{k}\right)\exp\left({x_1}^{k}+{x_2}^{k}+...+{x_n}^{k}\right) \frac{t^k}{k}
  \end{equation}
  with respect to $t$, then equating like powers of $t$, leads to the set of simultaneous equations for any positive integer $k$
\begin{equation}   \nonumber 
  1\times{_{n}}{B_{1}}
  = \left(1-a\right)\exp\left({x_1}+{x_2}+...+{x_n}\right),
\end{equation}
\begin{equation}   \nonumber
  2\times{_{n}}{B_{2}} = \left(1-a^{2}\right)\exp\left({x_1}^{2}+{x_2}^{2}+...+{x_n}^{2}\right) + {_{n}}{B_{1}} \times \left(1-a\right)\exp\left({x_1}+{x_2}+...+{x_n}\right),
\end{equation}
\begin{equation}  \nonumber
  3\times{_{n}}{B_{3}} = \left(1-a^{3}\right)\exp\left({x_1}^{3}+{x_2}^{3}+...+{x_n}^{3}\right) + {_{n}}{B_{1}} \times \left(1-a^{2}\right)\exp\left({x_1}^{2}+{x_2}^{2}+...+{x_n}^{2}\right)
\end{equation}
\begin{equation}  \nonumber
   + {_{n}}{B_{2}}\times \left(1-a\right)\exp\left({x_1}+{x_2}+...+{x_n}\right), \\
\end{equation}
\begin{equation}  \nonumber
  4\times{_{n}}{B_{4}} = \left(1-a^{4}\right)\exp\left({x_1}^{4}+{x_2}^{4}+...+{x_n}^{4}\right) + {_{n}}{B_{1}}\times \left(1-a^{3}\right)\exp\left({x_1}^{3}+{x_2}^{3}+...+{x_n}^{3}\right) \\
\end{equation}
\begin{equation}  \nonumber
  + {_{n}}{B_{2}}\times \left(1-a^{2}\right)\exp\left({x_1}^{2}+{x_2}^{2}+...+{x_n}^{2}\right) + {_{n}}{B_{3}} \times \left(1-a\right)\exp\left({x_1}+{x_2}+...+{x_n}\right), \\
\end{equation}
   etc., down to the $k$th equation
   \begin{equation}  \nonumber
  k\times{_{n}}{B_{k}}
\end{equation}
\begin{equation}   \nonumber
  = \left(1-a^{k}\right)\exp\left({x_1}^{k}+{x_2}^{k}+...+{x_n}^{k}\right)
\end{equation}
\begin{equation}  \nonumber
+ {_{n}}{B_{1}}\times \left(1-a^{k-1}\right)\exp\left({x_1}^{k-1}+{x_2}^{k-1}+...+{x_n}^{k-1}\right) + ... \\
\end{equation}
\begin{equation}  \nonumber
  ... + {_{n}}{B_{k-2}}\times \left(1-a^{2}\right)\exp\left({x_1}^{2}+{x_2}^{2}+...+{x_n}^{2}\right)
\end{equation}
\begin{equation}  \nonumber
+ {_{n}}{B_{k-1}}\times \left(1-a\right)\exp\left({x_1}+{x_2}+...+{x_n}\right). \\
\end{equation}
We note that for the above, the first equation has one unknown, the first and second equations together are \emph{two equations in two unknowns}, the first, second and third equations are \emph{three equations in three unknowns}, and so on. The set of $k$ equations here is a solvable set of \emph{k equations in k unknowns}, and is solved by Cramer's Rule involving determinants. Applying this rule yields the theorem.
\end{proof}

The above proof is essentially the same as that given earlier here, simply again applying Faà de Bruno's theorem.

We write theorem \ref{thm 4.1} as a single equation identity as follows.

\begin{equation}\label{5.1}
  G_n(x_1,x_2,x_3,...,x_n;a,t)= \prod_{\alpha_1,\alpha_2,...,\alpha_n \geq0} \left(\frac{1-x_1^{\alpha_1} x_2^{\alpha_2} ... x_n^{\alpha_n}at}{1-x_1^{\alpha_1} x_2^{\alpha_2} ... x_n^{\alpha_n}t}\right)^{\frac{1}{{\alpha_1}! {\alpha_2}! ... {\alpha_n}!}} \\
\end{equation}
\begin{equation}  \nonumber
  = 1 + \left(1-a\right)\exp\left({x_1}+{x_2}+...+{x_n}\right) \frac{t}{1!}
\end{equation}
\begin{equation}  \nonumber
+ \begin{vmatrix}
    \left(1-a\right)\exp\left({x_1}+{x_2}+...+{x_n}\right) & -1 \\
    \left(1-a^{2}\right)\exp\left({x_1}^{2}+{x_2}^{2}+...+{x_n}^{2}\right) & \left(1-a\right)\exp\left({x_1}+{x_2}+...+{x_n}\right) \\
  \end{vmatrix} \frac{t^2}{2!}
  \end{equation}
\begin{equation}  \nonumber
+ \begin{tiny}\begin{vmatrix}
    \left(1-a\right)\exp\left({x_1}+{x_2}+...+{x_n}\right) & -1 & 0 \\
    \left(1-a^{2}\right)\exp\left({x_1}^{2}+{x_2}^{2}+...+{x_n}^{2}\right) & \left(1-a\right)\exp\left({x_1}+{x_2}+...+{x_n}\right) & -2 \\
    \left(1-a^{3}\right)\exp\left({x_1}^{3}+{x_2}^{3}+...+{x_n}^{3}\right) & \left(1-a^{2}\right)\exp\left({x_1}^{2}+{x_2}^{2}+...+{x_n}^{2}\right) & \left(1-a\right)\exp\left({x_1}+{x_2}+...+{x_n}\right) \\
  \end{vmatrix} \frac{t^3}{3!}\end{tiny}
  \end{equation}
\begin{equation}  \nonumber
+ \begin{tiny}\begin{vmatrix}
    \left(1-a\right)\exp\left(\sum_{j=1}^{n}{x_j}\right) & -1 & 0 & 0 \\
    \left(1-a^{2}\right)\exp\left(\sum_{j=1}^{n}{x_j}^2\right) & \left(1-a\right)\exp\left(\sum_{j=1}^{n}{x_j}\right) & -2 & 0 \\
    \left(1-a^{3}\right)\exp\left(\sum_{j=1}^{n}{x_j}^3\right) & \left(1-a^{2}\right)\exp\left(\sum_{j=1}^{n}{x_j}^2\right) & \left(1-a\right)\exp\left(\sum_{j=1}^{n}{x_j}\right) & -3 \\
    \left(1-a^{4}\right)\exp\left(\sum_{j=1}^{n}{x_j}^4\right) & \left(1-a^{3}\right)\exp\left(\sum_{j=1}^{n}{x_j}^3\right) & \left(1-a^{2}\right)\exp\left(\sum_{j=1}^{n}{x_j}^2\right) & \left(1-a\right)\exp\left(\sum_{j=1}^{n}{x_j}\right) \\
  \end{vmatrix} \frac{t^4}{4!}\end{tiny}
  \end{equation}
\begin{equation}  \nonumber
 + etc.
\end{equation}

The following $1$-space example of equation (\ref{5.1}) is a variant on the $q$-binomial theorem itself. This result is not in the literature, and is equivalent to a statement about weighted integer partitions.
\begin{corollary}\label{5.2b}
\begin{equation}\label{5.2}
  G_1(q;a,t)= \prod_{k=0}^\infty \left(\frac{1-atq^k}{1-tq^k}\right)^{1/k!}
  = 1 + (1-a)e^q \frac{t}{1!}
    \end{equation}
  \begin{equation}  \nonumber
+ \begin{vmatrix}
    (1-a)e^q & -1 \\
    (1-a^2)e^{q^2} & (1-a)e^q \\
  \end{vmatrix} \frac{t^2}{2!}
+ \begin{vmatrix}
    (1-a)e^q & -1 & 0 \\
    (1-a^2)e^{q^2} & (1-a)e^q & -2 \\
    (1-a^3)e^{q^3} & (1-a^2)e^{q^2} & (1-a)e^q \\
  \end{vmatrix} \frac{t^3}{3!}
      \end{equation}
  \begin{equation}  \nonumber
+ \begin{vmatrix}
    (1-a)e^q & -1 & 0 & 0 \\
    (1-a^2)e^{q^2} & (1-a)e^q & -2 & 0 \\
    (1-a^3)e^{q^3} & (1-a^2)e^{q^2} & (1-a)e^q & -3 \\
    (1-a^4)e^{q^4} & (1-a^3)e^{q^3} & (1-a^2)e^{q^2} & (1-a)e^q \\
  \end{vmatrix} \frac{t^4}{4!}
+ etc.
\end{equation}
  \end{corollary}
    An interesting case of this is with $a=0$ so then

\begin{corollary}\label{5.3b}
\begin{equation}\label{5.3}
  G_1(q;0,t)= \prod_{k=0}^\infty \left(\frac{1}{1-tq^k}\right)^{1/k!}
   \end{equation}
  \begin{equation}  \nonumber
    = 1 + e^q \frac{t}{1!}
 + \begin{vmatrix}
    e^q & -1 \\
    e^{q^2} & e^q \\
  \end{vmatrix} \frac{t^2}{2!}
+ \begin{vmatrix}
    e^q & -1 & 0 \\
    e^{q^2} & e^q & -2 \\
    e^{q^3} & e^{q^2} & e^q \\
  \end{vmatrix} \frac{t^3}{3!}
 + \begin{vmatrix}
    e^q & -1 & 0 & 0 \\
    e^{q^2} & e^q & -2 & 0 \\
    e^{q^3} & e^{q^2} & e^q & -3 \\
    e^{q^4} & e^{q^3} & e^{q^2} & e^q \\
  \end{vmatrix} \frac{t^4}{4!}
+ etc.
\end{equation}
  \begin{equation}  \nonumber
    = 1 + e^q \frac{t}{1!}
 + \left(e^{q^2} + e^{2 q}\right) \frac{t^2}{2!}
+ \left( 2 e^{q^3} + 3 e^{q^2 + q} + e^{3q} \right) \frac{t^3}{3!}
\end{equation}
  \begin{equation}  \nonumber
 + \left( 6e^{q^4} + 8e^{q^3+q} + 3e^{2q^2} + 6e^{q^2+2q} + e^{4q}  \right) \frac{t^4}{4!}
+ etc.
\end{equation}
  \end{corollary}

Inferences from limiting cases of the coefficients lead us to affirm as checks and balances for example, that:

\begin{equation}  \label{5.4}
\begin{vmatrix}
    a & -1 & 0 & 0 \\
    a & a & -2 & 0 \\
    a & a & a & -3 \\
    a & a & a & a \\
  \end{vmatrix}
  = a(a+1)(a+2)(a+3), \\
\end{equation}

and

\begin{equation}  \label{5.5}
\begin{vmatrix}
    a & -1 & 0 & 0 \\
    a^2 & a & -2 & 0 \\
    a^3 & a^2 & a & -3 \\
    a^4 & a^3 & a^2 & a \\
  \end{vmatrix}
  = 4! a^4. \\
\end{equation}

Both (5.4) and (5.5) are easy to prove, and extend to generalized order determinants, and are consistent with application of Faà De Bruno's formula. It might be a worthwhile thing to publish a comprehensive list of cases of the types (5.4) and (5.5) applicable to the generating functions in our paper and other possible follow-up papers related to analysis of classes of vector partitions.

Next, repeating a variation on the recurrence logic from start of section 5, we see that

\begin{corollary} \label{5.6a}
\begin{equation}\label{5.6}
  G_1(q;0;t)/G_1(q^2;0;t^2) = \prod_{k=0}^\infty \left(1+tq^k\right)^{1/k!}
   \end{equation}
  \begin{equation}  \nonumber
    = 1 + e^q \frac{t}{1!}
 + \begin{vmatrix}
    -e^q & -1 \\
    e^{q^2} & -e^q \\
  \end{vmatrix} \frac{t^2}{2!}
+ \begin{vmatrix}
    -e^q & -1 & 0 \\
    e^{q^2} & -e^q & -2 \\
    -e^{q^3} & e^{q^2} & -e^q \\
  \end{vmatrix} \frac{t^3}{3!}
 + \begin{vmatrix}
    -e^q & -1 & 0 & 0 \\
    e^{q^2} & -e^q & -2 & 0 \\
    -e^{q^3} & e^{q^2} & -e^q & -3 \\
    e^{q^4} & -e^{q^3} & e^{q^2} & -e^q \\
  \end{vmatrix} \frac{t^4}{4!}
+ etc.
\end{equation}
  \begin{equation}  \nonumber
    = 1 + e^q \frac{t}{1!}
 + \left(e^{2q} - e^{q^2}\right) \frac{t^2}{2!}
+ \left( 2 e^{q^3} - 3 e^{q^2 + q} + e^{3 q} \right) \frac{t^3}{3!}
\end{equation}
  \begin{equation}  \nonumber
 + \left( e^{4q} + 3 e^{2q^2} - 6 e^{q^4} - 6 e^{2q + q^2} + 8 e^{q + q^3}  \right) \frac{t^4}{4!}
+ etc.
\end{equation}
  \end{corollary}

The statement of theorem \ref{thm 4.1} and it's form in equation (5.1) reminds us that considerations associated with
oscillating functions near the boundaries of convergence might apply or at least be worth covering off for exclusion
for such results. For an account of this phenomenon see Hardy and Littlewood’s paper \cite{gH1974c}
on the so-called “high indices theorem”. In particular, the reader should heed the early papers by
Hardy (see \cite{gH1974b} and \cite{gH1974a}) where it is clear that the behaviour of our functions near their poles
and the radii of convergence should be examined and
studied in order to eliminate the possibility that our formulas do not oscillate wildly in certain neighborhoods,
or are only picking up major terms and omitting complicated oscillating minor terms.

\begin{corollary}\label{5.7b}
The 2-space example of equation (\ref{5.1}) is
\begin{equation}\label{5.7}
  G_2(x,y;a,t)= \prod_{j,k\geq0}\left(\frac{1-x^j y^k at}{1-x^j y^k t}\right)^{\frac{1}{j! k!}} \\
\end{equation}
  \begin{equation}  \nonumber
= 1 + (1-a)e^{x+y} \frac{t}{1!} + \begin{vmatrix}
    (1-a)e^{x+y} & -1 \\
    (1-a^2)e^{x^2+y^2} & (1-a)e^{x+y} \\
  \end{vmatrix} \frac{t^2}{2!}
\end{equation}
  \begin{equation}  \nonumber
  + \begin{vmatrix}
    (1-a)e^{x+y} & -1 & 0 \\
    (1-a^2)e^{x^2+y^2} & (1-a)e^{x+y} & -2 \\
    (1-a^3)e^{x^3+y^3} & (1-a^2)e^{x^2+y^2} & (1-a)e^{x+y} \\
  \end{vmatrix} \frac{t^3}{3!}
      \end{equation}
  \begin{equation}  \nonumber
+ \begin{vmatrix}
    (1-a)e^{x+y} & -1 & 0 & 0 \\
    (1-a^2)e^{x^2+y^2} & (1-a)e^{x+y} & -2 & 0 \\
    (1-a^3)e^{x^3+y^3} & (1-a^2)e^{x^2+y^2} & (1-a)e^{x+y} & -3 \\
    (1-a^3)e^{x^4+y^4} & (1-a^3)e^{x^3+y^3} & (1-a^2)e^{x^2+y^2} & (1-a)e^{x+y} \\
  \end{vmatrix} \frac{t^4}{4!}
+ etc.
\end{equation}
  \end{corollary}

  The function $G_2(x,y;a,t)$ is interesting for a variety of reasons, and can lead to new and interesting combinatorial analysis involving weighted partitions of vectors.  In other words, we can examine stepping stone jumps between integer lattice points in 2-space, where each jump to the next stepping stone involves carrying a weight assigned to the coefficient, and our "vector partition science" will be based upon the total weight carried in order to combine all possible stepping stone jumps assigned by rule. In 3-space this same weighted stepping stone combinatorial analysis will also apply, and the 3-space version of (5.1) given now is:

\begin{corollary}\label{5.8b}
The 3-space example of equation (\ref{5.1}) is
\begin{equation}\label{5.8}
  G_3(x,y,z;a,t)= \prod_{h,j,k\geq0}\left(\frac{1-x^h y^j z^k at}{1-x^h y^j z^k t}\right)^{\frac{1}{h! j! k!}} \\
\end{equation}
  \begin{equation}  \nonumber
= 1 + + (1-a)e^{x+y+z} \frac{t}{1!} + \begin{vmatrix}
    (1-a)e^{x+y+z} & -1 \\
    (1-a^2)e^{x^2+y^2+z^2} & (1-a)e^{x+y+z} \\
  \end{vmatrix} \frac{t^2}{2!}
\end{equation}
  \begin{equation}  \nonumber
  + \begin{vmatrix}
    (1-a)e^{x+y+z} & -1 & 0 \\
    (1-a^2)e^{x^2+y^2+z^2} & (1-a)e^{x+y+z} & -2 \\
    (1-a^3)e^{x^3+y^3+z^3} & (1-a^2)e^{x^2+y^2+z^2} & (1-a)e^{x+y+z} \\
  \end{vmatrix} \frac{t^3}{3!}
      \end{equation}
  \begin{equation}  \nonumber
+ \begin{vmatrix}
    (1-a)e^{x+y+z} & -1 & 0 & 0 \\
    (1-a^2)e^{x^2+y^2+z^2} & (1-a)e^{x+y+z} & -2 & 0 \\
    (1-a^3)e^{x^3+y^3+z^3} & (1-a^2)e^{x^2+y^2+z^2} & (1-a)e^{x+y+z} & -3 \\
    (1-a^3)e^{x^4+y^4+z^4} & (1-a^3)e^{x^3+y^3+z^3} & (1-a^2)e^{x^2+y^2+z^2} & (1-a)e^{x+y+z} \\
  \end{vmatrix} \frac{t^4}{4!}
+ etc.
\end{equation}
  \end{corollary}

We reiterate here that the identities in the present paper follow the approach taken by Andrews \cite[chapter 2]{gA1976}, who develops his theory of integer partitions starting from an historical viewpoint and uses the $q$-binomial theorem to derive the Euler identities for the number of partitions of a positive integer. It would be interesting if (5.6) or (5.7) could yield higher Jacobi triple product analogies, or higher 2-space and 3-space type Heine $q$-series summations and transforms, and other 2-space and 3-space $q$-series sum-product identities relevant to vector partitions.

\section{An $n$-space $q$-binomial functional equation.} \label{S:Functional equation}

We now state a functional equation involving the $n$-space $q$-binomial theorem defined in \S2, highlighting the 2-space and 3-space cases as examples.
\begin{theorem} \label{6.1a}
  If the left side of (3.1) is for the moment considered only as a function, $F(t)$, of $t$, then
  \begin{equation}\label{6.1}
    \frac{1-at}{1-t} = \frac{F(t)\prod_{S_2}F(x_{i_1} x_{i_2} t) \prod_{S_4} F(x_{i_1} x_{i_2} x_{i_3} x_{i_4} t) ... }
    {\prod_{S_1}F(x_{i_1} t) \prod_{S_3} F(x_{i_1} x_{i_2} x_{i_3} t)\prod_{S_5}F(x_{i_1} x_{i_2} x_{i_3} x_{i_4} x_{i_5} t) ...  }
  \end{equation}
  the numerator products each being over the $n$th order even symmetric combinations whilst the denominator products are over the $n$th order odd symmetric combinations of the variables.  We have used the notation $S_1, S_2, S_3, S_4, ...$ to denote the respective positive integer lattice sets operating on the right side products here, such that
\begin{eqnarray}
 \nonumber 
  S_1 &=& 1 \leq i_1 \leq n; \\
 \nonumber
  S_2 &=& 1 \leq i_1 \leq n, 1 \leq i_2 \leq n; \\
 \nonumber
  S_3 &=& 1 \leq i_1 \leq n, 1 \leq i_2 \leq n, 1 \leq i_3 \leq n; \\
 \nonumber
  S_4 &=& 1 \leq i_1 \leq n, 1 \leq i_2 \leq n, 1 \leq i_3 \leq n, 1 \leq i_4 \leq n; \\
 \nonumber
   etc.,
\end{eqnarray}
and the right side of (\ref{4.1}) is a finite product of such functions going up to the  $n$th  order symmetric function product of terms.
\end{theorem}

It is worth noting (and spelling out) that we have used the abbreviation $F(t)=F_n(x_1,x_2,x_3,...,x_n;a,t)$
and the right side of (\ref{4.1}) is a finite product of such functions going up to the $n$th order symmetric function product of terms. It is also seen here that (\ref{4.1}) is a functional equation for a \emph{generalized right side of (\ref{1.1}}) form of the $q$-binomial product. As an illustration we next give the $n=2$ and $n=3$ examples of the functional equation (\ref{4.1}), with $x_1 \rightarrow x$, $x_2 \rightarrow y$, $x_3 \rightarrow z$. With these variables, we assert
\begin{corollary}
  If for $\left|x\right|<1, \left|y\right|<1, \left|z\right|<1$, and all complex numbers $a$ and $t$,
  \begin{eqnarray*}
    F_2(x,y;a,t) &=& \prod_{j,k\geq0} \frac{1-x^j y^k at}{1-x^j y^k t},  \\
    F_3(x,y,z;a,t) &=& \prod_{j,k,h\geq0} \frac{1-x^j y^k z^h at}{1-x^j y^k z^h t};
  \end{eqnarray*}
  then
  \begin{eqnarray*}
    \frac{1-at}{1-t} &=& \frac{F_2(x,y;a,t)F_2(x,y;a,xyt)}{F_2(x,y;a,xt)F_2(x,y;a,yt)},  \\
    \frac{1-at}{1-t} &=& \frac{F_3(x,y,z;a,t)F_3(x,y,z;a,xyt)F_3(x,y,z;a,xzt)F_3(x,y,z;a,yzt)}
    {F_3(x,y,z;a,xt)F_3(x,y,z;a,yt)F_3(x,y,z;a,zt)F_3(x,y,z;a,xyzt)}.
  \end{eqnarray*}
\end{corollary}

\section{A natural application to the Visible Point Vector (VPV) identities.} \label{S:VPV identities}

In the 1990s and up to 2000 the author published a series of papers introducing the so-called Visible Point Vector (VPV) identities. In these papers by the author (see for example Campbell \cite{gC1994a}, \cite{gC1994b}, \cite{gC1998} and \cite{gC2000}), the identities given in the present section were published. They attracted scant attention at the time, and were seen as curiosities. They are however, generating functions for weighted vector partitions, not too different conceptually to the ones presented in the sections so far in the current paper.

So, we establish the following conventions for the VPV identities.
\begin{definition}
  We use the notation $\left(x_1, x_2, x_3,...,x_n\right)$, to mean “the greatest common divisor of all of
$x_1, x_2, x_3,...,x_n$ together; the same as $\gcd\left(x_1, x_2, x_3,...,x_n\right)$”.
It is important to distinguish between this and the ordered $n$-tuple utilized for the vector
$\langle x_1, x_2, x_3,...,x_n \rangle$. In either case we will be concerned with lattice points
in the relevant Euclidean space, hence any vector or gcd will be over integer coordinates.
\end{definition}

\begin{definition}
  Any Euclidean vector $\langle x_1, x_2, x_3,...,x_n \rangle$ for which
$\left( x_1, x_2, x_3,...,x_n \right)=1$ we call a visible point vector, abbreviated VPV.
\end{definition}

\begin{theorem}   \label{7.1a}
  \textbf{The first hyperquadrant VPV identity.} If $i = 1, 2, 3,...,n$ then for each $x_i \in \mathbb{C}$ such that $|x_i|<1$ and $b_i \in \mathbb{C}$ such that $\sum_{i=1}^{n}b_i = 1$,
  \begin{equation}   \label{7.1}
    \prod_{\substack{ (a_1,a_2,...,a_n)=1 \\ a_1,a_2,...,a_n \geq 1}} \left( \frac{1}{1-{x_1}^{a_1}{x_2}^{a_2}{x_3}^{a_3}\cdots{x_n}^{a_n}} \right)^{\frac{1}{{a_1}^{b_1}{a_2}^{b_2}{a_3}^{b_3}\cdots{a_n}^{b_n}}}
  \end{equation}
  \begin{equation}  \nonumber
  = \exp\left\{ \prod_{i=1}^{n} \left( \sum_{j=1}^{\infty} \frac{{x_i}^j}{j^{b_i}} \right)\right\}.
  \end{equation}
\end{theorem}

There follow numerous example corollaries of this theorem, all of them susceptible to the combinatorial analysis of the previous sections. However, firstly we give the lemma and proof underpinning theorem 7.1.

\begin{lemma}
  Consider an infinite region raying out of the origin in any Euclidean
vector space. The set of all lattice point vectors apart from the origin in that region is
precisely the set of positive integer multiples of the VPVs in that region.
\end{lemma}
\begin{proof}
  Each VPV will have integer coordinates whose greatest common divisor
is unity. Viewed from the origin, all other lattice points are obscured behind the VPV end
points. If $x$ is a VPV in the region then all vectors in that region from the origin with direction
of $x$ preserved are enumerated by a sequence $1x, 2x, 3x,...$, and the greatest
common divisor of the components of $nx$ is clearly $n$. This is because if the scalar $n$ is
non-integer at least one of the coordinates of $nx$ would be a non-integer. Therefore, if
the VPVs in the region are countably given by $x_1, x_2, x_3, ...$, then all lattice point vectors
from the origin in the region are
$1x_1,2x_1,3x_1,...; 1x_2,2x_2,3x_2,...; 1x_3,2x_3,3x_3,...$ etc.
Completion of the proof comes with recognition that the set of all VPVs in any \emph{rayed
from the origin} region in any Euclidean vector space is a countable set. Proof of this last
assertion is by induction on the dimension, knowing that the lattice points are countable
in any two dimensional region. As we count each lattice point vector in the desired region
we decide whether it is a VPV simply by observing whether its coordinates are relatively
prime as a whole.
\end{proof}

This then brings us to the proof of theorem 7.1.
\begin{proof}
  We start with the multiple sum

  \begin{equation} \nonumber
    \sum_{ a_1,a_2,...,a_n \in \mathbb{Z}^+} \frac{{x_1}^{a_1}{x_2}^{a_2}{x_3}^{a_3}\cdots{x_n}^{a_n}}{{a_1}^{b_1}{a_2}^{b_2}{a_3}^{b_3}\cdots{a_n}^{b_n}}
    = \prod_{i=1}^{n} \sum_{j=1}^{\infty} \frac{{x_i}^j}{j^{b_i}}
    \end{equation}
  which, due to Lemma 7.1, also equals, letting $b = \sum_{i=1}^{n} b_i$,
 \begin{tiny} \begin{equation} \nonumber
    \sum_{\substack{ (a_1,a_2,...,a_n)=1 \\ a_1,a_2,...,a_n \geq 1}}
           \left( \frac{{x_1}^{a_1}{x_2}^{a_2}\cdots{x_n}^{a_n}}{1^b}
            + \frac{({x_1}^{a_1}{x_2}^{a_2}\cdots{x_n}^{a_n})^2}{2^b}
            + \frac{({x_1}^{a_1}{x_2}^{a_2}\cdots{x_n}^{a_n})^3}{3^b}
             + \cdots\right)
             \frac{1}{{a_1}^{b_1}{a_2}^{b_2}{a_3}^{b_3}\cdots{a_n}^{b_n}}
  \end{equation}\end{tiny}
  \begin{equation} \nonumber
    = \sum_{\substack{ (a_1,a_2,...,a_n)=1 \\ a_1,a_2,...,a_n \geq 1}}
             \frac{- \log( 1 - {x_1}^{a_1}{x_2}^{a_2}\cdots{x_n}^{a_n})}{{a_1}^{b_1}{a_2}^{b_2}{a_3}^{b_3}\cdots{a_n}^{b_n}}
  \end{equation}
  Exponentiating both sides then yields Theorem 7.1.
\end{proof}
The cases of theorem 7.1 with $n = 2, n = 3$, are stated easily in the forms,

\begin{corollary} If $|y|<1, |z|<1,$ and $s+t=1$, then
  \begin{equation}   \label{7.2}
    \prod_{\substack{ (a,b)=1 \\ a,b \geq 1}} \left( \frac{1}{1-y^a z^b} \right)^{\frac{1}{a^s b^t}}
    = \exp\left\{ \left( \sum_{i=1}^{\infty} \frac{y^i}{i^s}\right) \left( \sum_{j=1}^{\infty} \frac{z^j}{j^t}\right) \right\}.
  \end{equation}
\end{corollary}

\begin{corollary} If $|x|<1, |y|<1, |z|<1 $ and $s+t+u=1$, then
  \begin{equation}   \label{7.3}
    \prod_{\substack{ (a,b,c)=1 \\ a,b,c \geq 1}} \left( \frac{1}{1-x^a y^b z^c} \right)^{\frac{1}{a^s b^t c^u}}
    = \exp\left\{ \left( \sum_{i=1}^{\infty} \frac{x^i}{i^s}\right) \left( \sum_{j=1}^{\infty} \frac{y^j}{j^t}\right)
    \left( \sum_{k=1}^{\infty} \frac{z^k}{k^u}\right)\right\}.
  \end{equation}
\end{corollary}
The reader will recognise the polylogarithm occurring in the right sides of (7.2) and (7.3). The limiting values of the polylogarithms being Riemann zeta functions implies interesting new identities such as,

\begin{equation}   \label{7.4}
    \prod_{\substack{ (a,b)=1 \\ a,b \geq 1}} \left( \frac{1}{1-y^a z^b} \right)^{\frac{1}{\surd{(ab)}}}
    = \exp\left\{ \left( \sum_{i=1}^{\infty} \frac{y^i}{\surd{i}}\right) \left( \sum_{j=1}^{\infty} \frac{z^j}{\surd{j}}\right) \right\},
  \end{equation}
  \begin{equation}   \label{7.5}
    \prod_{\substack{ (a,b,c)=1 \\ a,b,c \geq 1}} \left( \frac{1}{1-x^a y^b z^c} \right)^{\frac{1}{{(abc)}^{\frac{1}{3}}}}
    = \exp\left\{ \left( \sum_{i=1}^{\infty} \frac{x^i}{i^{\frac{1}{3}}}\right) \left( \sum_{j=1}^{\infty} \frac{y^j}{j^{\frac{1}{3}}}\right)
    \left( \sum_{k=1}^{\infty} \frac{z^k}{k^{\frac{1}{3}}}\right)\right\}.
  \end{equation}

  We are reminded by (7.2) of the functional equation due originally to Riemann
in his famous paper \cite{gR1859} on the Riemann zeta function. Both have the $s+t=1$ caveat. Riemann's zeta function reflection
formula is equivalent to

\begin{equation}\label{7.6}
  \Gamma\left( \frac{s}{2} \right)\pi^{-s/2}\zeta(s) = \Gamma\left( \frac{t}{2} \right)\pi^{-t/2}\zeta(t)
\end{equation}
  where $s+t=1$, but equation (7.2) is quite a different relationship in a context amenable to the critical line Riemann zeta function $\zeta\left( \frac{1}{2} + i t\right)$ for nontrivial zeroes.

There are several further corollary cases that we can state here, that may be susceptible to the analysis of the earlier sections. There are natural and simple cases of Theorem 7.1 to consider. Let us first enlarge the theorem’s positive coordinate hyperquadrant to include lattice points on each axis except for the highest or $n$th dimension. In other words, the product operator for variable $z$ on each left side of (7.7) to (7.10) runs over each integer 1, 2, 3,... whereas for the non$z$ variables $v, w, x, y$, the product is over 0, 1, 2, 3,....  Thus we can easily obtain the following infinite products involving VPVs in the combinatorial interpretations.

\begin{corollary} For each of $|v|, |w|, |x|, |y|, |z|<1,$
  \begin{equation}\label{7.7}
    \prod_{\substack{(a,b)=1 \\ a\geq0,b>0}} \left( 1-y^a z^b \right)^{\frac{1}{b}}
    = (1-z)^{\frac{1}{1-y}},
  \end{equation}
    \begin{equation}\label{7.8}
    \prod_{\substack{(a,b,c)=1 \\ a,b\geq0,c>0}} \left( 1-x^a y^b z^c \right)^{\frac{1}{c}}
    = (1-z)^{\frac{1}{(1-x)(1-y)}},
  \end{equation}
    \begin{equation}\label{7.9}
    \prod_{\substack{(a,b,c,d)=1 \\ a,b,c\geq0,d>0}} \left( 1-w^a x^b y^c z^d \right)^{\frac{1}{d}}
    = (1-z)^{\frac{1}{(1-w)(1-x)(1-y)}},
  \end{equation}
    \begin{equation}\label{7.10}
    \prod_{\substack{(a,b,c,d,e)=1 \\ a,b,c,d\geq0,e>0}} \left( 1-v^a w^b x^c y^d z^e \right)^{\frac{1}{e}}
    = (1-z)^{\frac{1}{(1-v)(1-w)(1-x)(1-y)}}.
  \end{equation}
  \end{corollary}

The above four infinite products and their reciprocals are worth deeper analysis as simple examples of weighted VPV partitions, reminiscent of the integer partition theorems. (\ref{7.7}) to (\ref{7.10}) are interesting to examine in the "near bijection" context that has been applied to the classical Euler pentagonal number theorem. This is a large topic probably beyond the scope of our present paper.

Let is take the example of equation (\ref{7.8}). The right side product is a case of the binomial theorem, which when applied gives us,

  \begin{equation}\label{7.11}
    \prod_{\substack{(a,b,c)=1 \\ a,b\geq0,c>0}} \left( 1-x^a y^b z^c \right)^{\frac{1}{c}}
   = (1-z)^{\frac{1}{(1-x)(1-y)}}
      \end{equation}
  \begin{equation}  \nonumber
   = 1 + \frac{(1-y)}{1!} \left(\frac{z}{1-x}\right)^1
    + \frac{(1-y)(x-y)}{2!} \left(\frac{z}{1-x}\right)^2
     \end{equation}
  \begin{equation}  \nonumber
     + \frac{(1-y)(x-y)(2x-y-1)}{3!} \left(\frac{z}{1-x}\right)^3
  \end{equation}
  \begin{equation}  \nonumber
    + \frac{(1-y)(x-y)(2x-y-1)(3x-y-2)}{4!} \left(\frac{z}{1-x}\right)^4 + \cdots
  \end{equation}

Looking closer at this, one sees that (\ref{7.11}) encodes a theorem about weighted VPVs, pertaining to visible from the origin points. For partitions of these VPVs in the first hyperquadrant of Euclidean 3-space, each vector $\langle a,b,c \rangle$ has integer coordinates that satisfy $a,b\geq0,c>0$. By a weighted partition, we mean a "\textit{\textbf{stepping stone jump while carrying a weight determined by a coefficient }}" from one integer lattice point to the next, jumping always "\textbf{\textit{away from the origin}}", that origin being the point $\langle 0,0,0 \rangle$. ie. The distance $\sqrt{a^2+b^2+c^2}$ from $\langle 0,0,0 \rangle$ to the starting point $\langle a,b,c \rangle$ of the jump is less than the distance $\sqrt{h^2+j^2+k^2}$ from $\langle 0,0,0 \rangle$ to the destination point $\langle h,j,k \rangle$ of the jump.

We here also give examples from the hyperpyramid VPV theorem given by the author in \cite{gC2000}, as the application of the determinant coefficient technique of our current paper is strikingly applicable and bearing some semblance to the $q$-binomial variants. Note that for each of (7.11) to (7.15) the left side products are taken over a set of integer lattice points inside an inverted hyperpyramid on the Euclidean cartesian space.

\begin{corollary} For $|y|, |z|<1,$
  \begin{equation}\label{7.12}
    \prod_{\substack{(a,b)=1 \\ a<b \\ a\geq0,b>0}} \left( \frac{1}{1-y^a z^b} \right)^{\frac{1}{b}}
    = \left(\frac{1-yz}{1-z}\right)^{\frac{1}{1-y}}
  \end{equation}
        \begin{equation}  \nonumber
= 1 + \frac{z}{1!} + \begin{vmatrix}
    1 & -1 \\
    \frac{1-y^2}{1-y} & 1 \\
  \end{vmatrix} \frac{z^2}{2!}
  + \begin{vmatrix}
    1 & -1 & 0 \\
    \frac{1-y^2}{1-y} & 1 & -2 \\
    \frac{1-y^3}{1-y} & \frac{1-y^2}{1-y} & 1 \\
  \end{vmatrix} \frac{z^3}{3!}
+ \begin{vmatrix}
    1 & -1 & 0 & 0 \\
    \frac{1-y^2}{1-y} & 1 & -2 & 0 \\
    \frac{1-y^3}{1-y} & \frac{1-y^2}{1-y} & 1 & -3 \\
    \frac{1-y^4}{1-y} & \frac{1-y^3}{1-y} & \frac{1-y^2}{1-y} & 1 \\
  \end{vmatrix} \frac{z^4}{4!}
+ etc.
\end{equation}
   \end{corollary}
  In this case it is fairly easy to find the Taylor coefficients for the (\ref{7.12}) right side function. Hence we get a closed form evaluation of the determinant coefficients. In Mathematica, and WolframAlpha one easily sees that the Taylor series is

  \begin{equation}  \nonumber
    \left(\frac{1-yz}{1-z}\right)^{\frac{1}{1-y}} = 1 + z + (y + 2) \frac{z^2}{2!} + (2 y^2 + 5 y + 6) \frac{z^3}{3!} + (6 y^3 + 17 y^2 + 26 y + 24) \frac{z^4}{4!}
  \end{equation}
  \begin{equation} \nonumber
   + (24 y^4 + 74 y^3 + 129 y^2 + 154 y + 120) \frac{z^5}{5!} + O(z^6)
  \end{equation}
  and that the expansion is encapsulated by
   $\sum_{n=0}^{\infty} c_n z^n$ where $c_0 = 1$, $c_1 = 1$ with the recurrence
   \begin{equation} \nonumber
   ny c_n + (n+2) c_{n+2} = (2 + n + y + ny) c_{n+1}.
   \end{equation}

   Incidentally, also in Mathematica, and WolframAlpha one easily sees that the code
   \begin{equation} \nonumber
   Det[\{1,-1,0,0\},\{(1-y^2)/(1-y),1,-2,0\},\{(1-y^3)/(1-y),(1-y^2)/(1-y),1,-3\}.
  \end{equation}
  \begin{equation} \nonumber
  \{(1-y^4)/(1-y),(1-y^3)/(1-y),(1-y^2)/(1-y),1\}]
   \end{equation}
   nicely verifies the coefficient given by

  \begin{equation}  \nonumber
\begin{vmatrix}
    1 & -1 & 0 & 0 \\
    \frac{1-y^2}{1-y} & 1 & -2 & 0 \\
    \frac{1-y^3}{1-y} & \frac{1-y^2}{1-y} & 1 & -3 \\
    \frac{1-y^4}{1-y} & \frac{1-y^3}{1-y} & \frac{1-y^2}{1-y} & 1 \\
  \end{vmatrix}
  = 6 y^3 + 17 y^2 + 26 y + 24. \\
\end{equation}

\begin{corollary} For each of $|x|, |y|, |z|<1,$
    \begin{equation}\label{7.13}
    \prod_{\substack{(a,b,c)=1 \\ a,b<c \\ a,b\geq0,c>0}} \left( \frac{1}{1-x^a y^b z^c} \right)^{\frac{1}{c}}
    = \left(\frac{(1-xz)(1-yz)}{(1-z)(1-xyz)}\right)^{\frac{1}{(1-x)(1-y)}}
  \end{equation}
      \begin{equation}  \nonumber
= 1 + \frac{z}{1!} + \begin{vmatrix}
    1 & -1 \\
    \frac{(1-x^2)(1-y^2)}{(1-x)(1-y)} & 1 \\
  \end{vmatrix} \frac{z^2}{2!}
  + \begin{vmatrix}
    1 & -1 & 0 \\
    \frac{(1-x^2)(1-y^2)}{(1-x)(1-y)} & 1 & -2 \\
    \frac{(1-x^3)(1-y^3)}{(1-x)(1-y)} & \frac{(1-x^2)(1-y^2)}{(1-x)(1-y)} & 1 \\
  \end{vmatrix} \frac{z^3}{3!}
      \end{equation}
  \begin{equation}  \nonumber
+ \begin{vmatrix}
    1 & -1 & 0 & 0 \\
    \frac{(1-x^2)(1-y^2)}{(1-x)(1-y)} & 1 & -2 & 0 \\
    \frac{(1-x^3)(1-y^3)}{(1-x)(1-y)} & \frac{(1-x^2)(1-y^2)}{(1-x)(1-y)} & 1 & -3 \\
    \frac{(1-x^4)(1-y^4)}{(1-x)(1-y)} & \frac{(1-x^3)(1-y^3)}{(1-x)(1-y)} & \frac{(1-x^2)(1-y^2)}{(1-x)(1-y)} & 1 \\
  \end{vmatrix} \frac{z^4}{4!}
+ etc.
\end{equation}
   \end{corollary}

\begin{corollary} For each of $|w|, |x|, |y|, |z|<1,$
    \begin{equation}\label{7.14}
    \prod_{\substack{(a,b,c,d)=1 \\ a,b,c<d \\ a,b,c\geq0,d>0}} \left( \frac{1}{1-w^a x^b y^c z^d} \right)^{\frac{1}{d}}
    = \left(\frac{(1-wz)(1-xz)(1-yz)(1-wxyz)}{(1-z)(1-wxz)(1-wyz)(1-xyz)}\right)^{\frac{1}{(1-w)(1-x)(1-y)}},
  \end{equation}
    \begin{equation}  \nonumber
= 1 + \frac{z}{1!} + \begin{vmatrix}
    1 & -1 \\
    \frac{(1-w^2)(1-x^2)(1-y^2)}{(1-w)(1-y)(1-z)} & 1 \\
  \end{vmatrix} \frac{z^2}{2!}
\end{equation}
  \begin{equation}  \nonumber
  + \begin{vmatrix}
    1 & -1 & 0 \\
    \frac{(1-w^2)(1-x^2)(1-y^2)}{(1-w)(1-y)(1-z)} & 1 & -2 \\
    \frac{(1-w^3)(1-x^3)(1-y^3)}{(1-w)(1-y)(1-z)} & \frac{(1-w^2)(1-x^2)(1-y^2)}{(1-w)(1-y)(1-z)} & 1 \\
  \end{vmatrix} \frac{z^3}{3!}
      \end{equation}
  \begin{equation}  \nonumber
+ \begin{vmatrix}
    1 & -1 & 0 & 0 \\
    \frac{(1-w^2)(1-x^2)(1-y^2)}{(1-w)(1-x)(1-y)} & 1 & -2 & 0 \\
    \frac{(1-w^3)(1-x^3)(1-y^3)}{(1-w)(1-x)(1-y)} & \frac{(1-w^2)(1-x^2)(1-y^2)}{(1-w)(1-x)(1-y)} & 1 & -3 \\
    \frac{(1-w^4)(1-x^4)(1-y^4)}{(1-w)(1-x)(1-y)} & \frac{(1-w^3)(1-x^3)(1-y^3)}{(1-w)(1-x)(1-y)} & \frac{(1-w^2)(1-x^2)(1-y^2)}{(1-w)(1-x)(1-y)} & 1 \\
  \end{vmatrix} \frac{z^4}{4!}
+ etc.
\end{equation}
   \end{corollary}

\begin{corollary} For each of $|v|, |w|, |x|, |y|, |z|<1,$
    \begin{equation}\label{7.15}
    \prod_{\substack{(a,b,c,d,e)=1 \\ a,b,c,d<e \\ a,b,c,d\geq0,e>0}} \left( \frac{1}{1-v^a w^b x^c y^d z^e} \right)^{\frac{1}{e}}
   \end{equation}
   \begin{gather}\nonumber
     =       \left(\frac{(1-vz)(1-wz)(1-xz)(1-yz)}{(1-z)(1-vwz)(1-vxz)(1-vyz)}\right)^{\frac{1}{(1-v)(1-w)(1-x)(1-y)}} \\  \nonumber
      \times \left(\frac{(1-vwxz)(1-vwyz)(1-vxyz)(1-wxyz)}{(1-wxz)(1-wyz)(1-xyz)(1-vwxyz)}\right)^{\frac{1}{(1-v)(1-w)(1-x)(1-y)}}.
   \end{gather}
    \begin{equation}  \nonumber
= 1 + \frac{z}{1!} + \begin{vmatrix}
    1 & -1 \\
    \frac{(1-v^2)(1-w^2)(1-x^2)(1-y^2)}{(1-v)(1-w)(1-x)(1-y)} & 1 \\
  \end{vmatrix} \frac{z^2}{2!}
\end{equation}
  \begin{equation}  \nonumber
  + \begin{vmatrix}
    1 & -1 & 0 \\
    \frac{(1-v^2)(1-w^2)(1-x^2)(1-y^2)}{(1-v)(1-w)(1-x)(1-y)} & 1 & -2 \\
    \frac{(1-v^3)(1-w^3)(1-x^3)(1-y^3)}{(1-v)(1-w)(1-x)(1-y)} & \frac{(1-v^2)(1-w^2)(1-x^2)(1-y^2)}{(1-v)(1-w)(1-x)(1-y)} & 1 \\
  \end{vmatrix} \frac{z^3}{3!}
      \end{equation}
  \begin{equation}  \nonumber
+ \begin{vmatrix}
    1 & -1 & 0 & 0 \\
    \frac{(1-v^2)(1-w^2)(1-x^2)(1-y^2)}{(1-v)(1-w)(1-x)(1-y)} & 1 & -2 & 0 \\
    \frac{(1-v^3)(1-w^3)(1-x^3)(1-y^3)}{(1-v)(1-w)(1-x)(1-y)} & \frac{(1-v^2)(1-w^2)(1-x^2)(1-y^2)}{(1-v)(1-w)(1-x)(1-y)} & 1 & -3 \\
    \frac{)1-v^4)(1-w^4)(1-x^4)(1-y^4)}{(1-v)(1-w)(1-x)(1-y)} & \frac{(1-v^3)(1-w^3)(1-x^3)(1-y^3)}{(1-v)(1-w)(1-x)(1-y)} & \frac{(1-v^2)(1-w^2)(1-x^2)(1-y^2)}{(1-v)(1-w)(1-x)(1-y)} & 1 \\
  \end{vmatrix} \frac{z^4}{4!}
+ etc.
 \end{equation}
    \end{corollary}

\section{A further example of the $n$-space variation of extended $q$-binomial theorem.} \label{S:section 8}

It is clear that we may apply the method to other vector partition generating functions.
An example is now given. The following theorem is based around the ideas associated with the elementary identity

\begin{equation}   \nonumber
  \left( 1 + x \right)\left( 1 + x^2 \right)\left( 1 + x^4 \right)\left( 1 + x^8 \right)\cdots = \frac{1}{1-x}.
\end{equation}

In fact the combinatorial interpretation of this is "\textit{\textbf{each positive integer is uniquely represented by a sum of distinct powers of 2}}".
So, we are here looking at an extension of this result in the

\begin{theorem}
 \begin{equation}\label{8.1}
    \prod_{k \geq 0} \left( \frac{1}{1- q t^{2^k}} \right)
    = \prod_{\substack{j,k \geq 0 \\ j \leq k}} \left( 1+ q^{2^j} t^{2^k} \right)
    = 1 + \sum_{k=1}^{\infty} A_k t^k
   \end{equation}
   where
   \begin{equation}  \nonumber
A_k = \frac{1}{k!}
   \begin{tiny}
   \begin{vmatrix}
    q & -1 & 0 & 0 & \cdots & 0 \\
    q^2+2q & q & -2 & 0 & \cdots & 0 \\
    q^3 & q^2+q & q & -3 & \cdots & 0 \\
    q^4+2q^2+4q & q^3 & q^2+q & q & \ddots & \vdots \\
    \vdots & \vdots & \vdots & \vdots & \ddots & -(k-1) \\
    \sum_{\substack{2^j|k \\ j \geq 0}} 2^j q^{k/2^j} & \sum_{\substack{2^j|(k-1) \\ j \geq 0}} 2^j q^{{(k-1)}/2^j} & \sum_{\substack{2^j|(k-2) \\ j \geq 0}} 2^j q^{{(k-2)}/2^j} & \sum_{\substack{2^j|(k-3) \\ j \geq 0}} 2^j q^{{(k-3)}/2^j} & \cdots & q \\
  \end{vmatrix} \end{tiny}
 \end{equation}
\end{theorem}

The combinatorial interpretation of (\ref{8.1}) is

\begin{theorem}
  If $B(j,k)$ is the number of vector partitions of $\langle j, k \rangle$ into distinct parts of kind
$\langle 2^a , 2^b \rangle$ in which $a \leq b$ with non-negative integers $a$ and $b$, then $B(j,k)$ equals also the number of partitions into
“unrestricted” parts of kind $\langle 1, 2^b \rangle$ in which $b$ is a non negative integer, and $B(j,k)$ is
the coefficient of $q^j t^k$ in (\ref{8.1}).
\end{theorem}

Each side of (\ref{8.1}) satisfies the equation $f(t)(1 - qt) = f(t^2)$ and this equation also
leads to a set of recurrences solvable using Cramer’s rule.

In Mathematica or Wolframalpha online we can easily check that

\begin{equation} \label{8.2}
  \prod_{k \geq 0} \left( \frac{1}{1- q t^{2^k}} \right)
  = 1 + q t + q (q + 1) t^2 + q^2 (q + 1) t^3 + q (q^3 + q^2 + q + 1) t^4
  \end{equation}
  \begin{equation}  \nonumber
  + q^2 (q^3 + q^2 + q + 1) t^5 + q^2 (q^4 + q^3 + q^2 + 2 q + 1) t^6
  \end{equation}
  \begin{equation}  \nonumber
  + q^3 (q^4 + q^3 + q^2 + 2 q + 1) t^7
  \end{equation}
  \begin{equation}  \nonumber
  + q (q^7 + q^6 + q^5 + 2 q^4 + 2 q^3 + q^2 + q + 1) t^8
\end{equation}
  \begin{equation}  \nonumber
  + q^2 (q^7 + q^6 + q^5 + 2 q^4 + 2 q^3 + q^2 + q + 1) t^9
\end{equation}
  \begin{equation}  \nonumber
  + q^2 (q^8 + q^7 + q^6 + 2 q^5 + 2 q^4 + 2 q^3 + 2 q^2 + 2 q + 1) t^{10} + \ldots
\end{equation}

Also, as a matter of interest, utilizing a form $\prod_{k=0}^\infty \prod_{j=0}^k (1+ q^{2^j} t^{2^k})$,
in Mathematica or Wolframalpha, the two product expansions in (\ref{8.1}) can be easily verified; both of them
yielding the series given in (\ref{8.2}). Therefore, to illustrate theorem 8.2, we give an arbitrary case for
the 2-space vector $\langle 3, 6 \rangle$:

\begin{corollary}
 $B(3,6)=2$ is the number of vector partitions of $\langle 3, 6 \rangle$ into distinct parts of kind
$\langle 2^a , 2^b \rangle$ in which $a \leq b$ with non-negative integers $a$ and $b$. The two partitions are $\langle 2, 4 \rangle + \langle 1, 2 \rangle$ and  $\langle 1, 4 \rangle + \langle 2, 2 \rangle$. Also $B(3,6)=2$ equals the number of partitions into
“unrestricted” parts of kind $\langle 1, 2^b \rangle$ in which $b$ is a non negative integer. The two partitions are $\langle 1, 2 \rangle + \langle 1, 2 \rangle + \langle 1, 2 \rangle$ and  $\langle 1, 4 \rangle + \langle 1, 1 \rangle + \langle 1, 1 \rangle$. Then also $B(3,6)=2$ is
the coefficient of $q^3 t^6$ in (\ref{8.1}).
\end{corollary}
The vector partitions $B(j,k)$ defined for theorem (\ref{8.1}) are easily visualized by the 2-space extension of a Ferrars graph of figure 1.

\begin{figure}
\centering
    \includegraphics[width=9cm,angle=0,height=9cm]{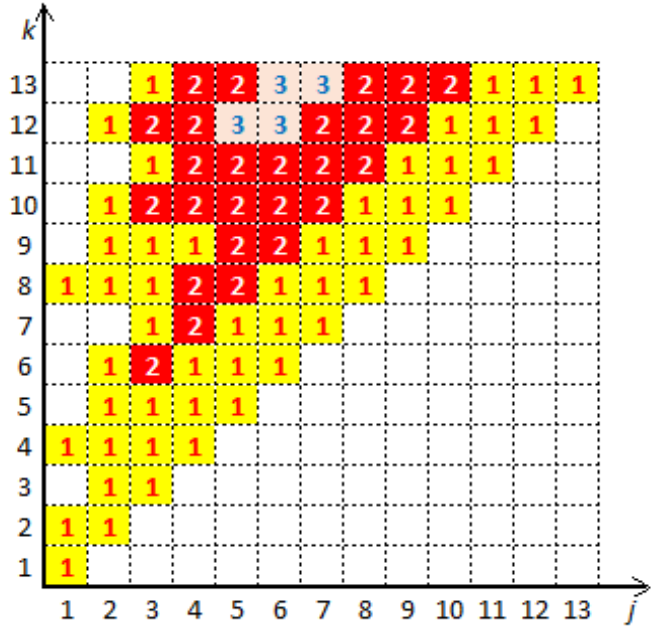}
  \caption{Number of vector partitions $B(j,k)$ encoded from theorem 8.2.}\label{Fig1}
\end{figure}

By now, the reader may surmise that the methods of this paper have set up the platform for possible interesting papers
delving further into the topic of vector partitions, and on other $n$-space multivariate infinite products.


\begin{thebibliography}{99}
\bibitem{mA1972}
ABRAMOWITZ, M., and STEGUN, I. Handbook of Mathematical Functions,
Dover Publications Inc., New York, 1972.
\bibitem{gA1976}
ANDREWS, G.E.  The Theory of Partitions,Addison-Wesley Publishing Company,
Advanced Book Program, Reading, Massachusetts, 1976.
\bibitem{tA1976}
APOSTOL, T.  Introduction to Analytic Number Theory,
Springer-Verlag, New York, 1976.
\bibitem{rB1982}
BAXTER, R. J. Exactly Solved Models in Statistical Mechanics, Academic Press,
New York, 1982.
\bibitem{gB1977}
BIRKHOFF, G. and MACLAINE, S. A survey of modern algebra, fourth ed., N.Y.,
Macmillan, 1977.
\bibitem{gC1994a}
CAMPBELL, G. B. “A new class of infinite products, and Euler's totient,” International Journal of Mathematics and Mathematical Sciences, vol. 17, no. 3, pp. 417-422, 1994. https://doi.org/10.1155/S0161171294000591.
\bibitem{gC2000}
CAMPBELL, G. B. “Infinite products over hyperpyramid lattices,” International Journal of Mathematics and Mathematical Sciences, vol. 23, no. 4, pp. 271-277, 2000. https://doi.org/10.1155/S0161171200000764.
\bibitem{gC1994b}
CAMPBELL, G. B. “Infinite products over visible lattice points,” International Journal of Mathematics and Mathematical Sciences, vol. 17, no. 4, pp. 637-654, 1994. https://doi.org/10.1155/S0161171294000918.
\bibitem{gC1998}
CAMPBELL, G. B. "A closer look at some new identities," International Journal of Mathematics and Mathematical Sciences, vol. 21, no. 3, pp. 581-586, 1998. https://doi.org/10.1155/S0161171298000805.
\bibitem{aC1893}
CAUCHY, A. M\'{e}moire sur les fonctions dont plusieurs \ldots , C. R. Acad. Sci. Paris,
T. XVII, p. 523, Oeuvres de Cauchy, 1re s\'{e}rie, T. VIII, Gauthier-Villars, Paris, 1893, 42-
50.
\bibitem{mC1966}
CHEEMA, M. S., Vector partitions and combinatorial identities, Math. Comp. 18,
1966 414-420.
\bibitem{mC1971}
CHEEMA, M. S. and MOTZKIN, T. S., Multipartitions and multipermutations, Proc.
Symp. Pure Math. 19, 1971, 37-39.
\bibitem{gG1990}
GASPER, G.  and  RAHMAN, M.  Basic Hypergeometric Series,
Encyclopedia of Mathematics and its Applications, Vol 35,
Cambridge University Press, (Cambridge - New York - Port Chester -
Melbourne - Sydney), 1990.
\bibitem{cG1813}
GAUSS, C.F.  Disquisitiones generales circa seriem infinitam
\ldots, Comm. soc. reg. sci. G\"{o}tt. rec., Vol II; reprinted in
Werke 3 (1876), pp. 123--162.
\bibitem{bG1963}
GORDON, B. Two theorems on multipartite partitions, J. London Math. Soc. 38,
1963, 459-464.
\bibitem{gH1974a}
HARDY, G. H. An extension of a theorem on oscillating series, Collected Papers, Vol
VI, Clarendon Press, Oxford, 1974, 500-506.
\bibitem{gH1974b}
HARDY, G. H. On certain oscillating series, Collected Papers, Vol VI, Clarendon
Press, Oxford, 1974, 146-167.
\bibitem{gH1974c}
HARDY, G. H., and LITTLEWOOD, J. E. A further note on the converse of Abel's
theorem. Collected Papers of Hardy, Vol VI, Clarendon Press, Oxford, 1974, 699-716.
\bibitem{eH1847}
HEINE, E.  Untersuchungen uber die Reihe ... , J. Reine angew.
Math. 34, 1847, 285-328.
\bibitem{eH1878}
HEINE, E.  Handbuch der Kugelfunctionen, Theorie und Andwendungen,
Vol. 1, Reimer, Berlin, 1878.
\bibitem{iM1995}
MACDONALD, I. G. Symmetric Functions And Hall Polynomials, 2nd ed., Oxford :
Clarendon Press ; New York : Oxford University Press, 1995.
\bibitem{gR1859}
RIEMANN, G. F. B. "\"{U}ber die Anzahl der Primzahlen unter einer gegebenen
Gr\"{o}sse." Monatsber. Königl. Preuss. Akad. Wiss. Berlin, 671-680, Nov. 1859.
\bibitem{oeiseuler}
SLOANE, N. J. A., The On-Line Encyclopedia of Integer Sequences (OEIS) Euler transform. ${\rm https://oeis.org/wiki/Euler\_transform}$.
\bibitem{nS2015}
SLOANE, N. J. A., The On-Line Encyclopedia of Integer Sequences (OEIS) sequence A061159 Numerators in expansion of Euler transform of b(n)=1/2 https://oeis.org/A061159.
\bibitem{nS2015a}
SLOANE, N. J. A., The On-Line Encyclopedia of Integer Sequences (OEIS) sequence A061160 Numerators in expansion of Euler transform of b(n)=1/3 https://oeis.org/A061160.
\bibitem{eW1956}
WRIGHT, E. M. Partitions of multipartite numbers, Proc. Amer. Math. Soc. 28, 1956,
880-890.
\end{thebibliography}
\end{document}